\documentclass[12pt]{amsart}
\usepackage{amsfonts}
\usepackage{amssymb}

\newtheorem{theorem}{Theorem}
\newtheorem{lemma}{Lemma}

\newtheorem{corollary}{Corollary}
\theoremstyle{definition}

\theoremstyle{remark}
\newtheorem*{remark}{Remark}

\DeclareMathOperator{\diam}{diam}\DeclareMathOperator{\conv}{conv}
\DeclareMathOperator{\Fix}{Fix}

\begin{document}
\title[On the fixed point property in direct sums]{On the fixed point property in direct sums of Banach spaces with
strictly monotone norms}
\author{Stanis{\l}aw Prus}
\author{Andrzej Wi\'{s}nicki}
\address{Stanis{\l}aw Prus, Institute of Mathematics, Maria Curie-Sk{\l}odowska University,
20-031 Lublin, Poland}
\email{bsprus@golem.umcs.lublin.pl}
\address{Andrzej Wi\'{s}nicki, Institute of Mathematics, Maria Curie-Sk{\l}odowska University,
20-031 Lublin, Poland} \email{awisnic@golem.umcs.lublin.pl}

\begin{abstract}
It is shown that if a Banach space $X$ has the weak Banach--Saks property and the weak fixed point
property for nonexpansive mappings and $Y$ satisfies property asymptotic (P) (which is weaker than
the condition $\text{\it WCS}\,( Y) >1$), then $X\oplus Y$ endowed with a strictly monotone norm
enjoys the weak fixed point property. The same conclusion is valid if $X$ admits a
$1$-unconditional basis.
\end{abstract}

\date{}
\maketitle

\section{Introduction}

One of the classic problems of metric fixed point theory concerns existence of fixed points of
nonexpansive mappings. Let $C$ be a nonempty bounded closed and convex subset of a Banach space
$X$. A mapping $T:C\rightarrow C$ is nonexpansive if
\begin{equation*}
\| Tx-Ty\| \leq \| x-y\|
\end{equation*}%
for all $x,y\in C$. A Banach space $X$ is said to have the fixed point property (FPP) if every such
mapping has a fixed point. Adding the assumption that $C$ is weakly compact in this condition, we
obtain the definition of the weak fixed point property (WFPP).

In 1965, F. Browder \cite{Br1} proved that Hilbert spaces have FPP. In the same year, Browder
\cite{Br2} and D. G\"{o}hde \cite{God} showed independently that uniformly convex spaces have FPP,
and W. A. Kirk \cite{Ki} proved a more general result stating that all Banach spaces with weak
normal structure have WFPP. Recall that a Banach space $X$ has weak normal structure if $r( C)
<\diam C$ for all weakly compact convex subsets $C$ of $X$ consisting of more than one point, where
$r( C)=\inf_{x\in C}\sup_{x\in C}\| x-y\|$ is the Chebyshev radius of $C$. There have been numerous
discoveries since then. In 1981, D. Alspach \cite{Al} showed an example of a nonexpansive
self-mapping defined on a weakly compact convex subset of $L_{1}[0,1]$ without a fixed point, and
B. Maurey \cite{Ma} used the Banach space ultraproduct construction to prove FPP for all reflexive
subspaces of $L_{1}[0,1]$ as well as WFPP for $c_{0}$, see also \cite{ELOS}. Maurey's method was
applied by P.-K. Lin \cite{Li1} who proved that every Banach space with a $1$-unconditional basis
enjoys WFPP. In 1997, P. Dowling and C. Lennard \cite{DoLe} proved that every nonreflexive subspace
of $L_{1}[0,1]$ fails FPP and they developed their techniques in the series of papers. For a fuller
discussion of metric fixed point theory we refer the reader to \cite{AyDoLo, GoKi, KiSi, KhKi}.

Major progress in fixed point problems for nonexpansive mappings has been made recently. In 2003
(published in 2006), J. Garc\'{\i}a Falset, E. Llor\'{e}ns Fuster and E. Mazcu\~{n}an Navarro
\cite{GLM2}, (see also \cite{MN}), solved a long-standing problem in the theory by proving FPP for
all uniformly nonsquare Banach spaces. In 2004, Dowling, Lennard and Turett \cite{DoLeTu} proved
that a nonempty closed bounded convex subset of $c_{0}$ has FPP if and only if it is weakly
compact. In a recent paper \cite{Li2}, P.-K. Lin showed that a certain renorming of $\ell _{1}$
enjoys FPP, thus solving another long-standing problem (FPP does not imply reflexivity).

The problem of whether FPP or WFPP is preserved under direct sums of Banach spaces has been
thoroughly studied since the 1968 Belluce--Kirk--Steiner theorem \cite{BKS}, which states that a
direct sum of two Banach spaces with normal structure, endowed with the maximum norm, also has
normal structure. In 1984, T. Landes \cite{La1} showed that normal structure is preserved under a
large class of direct sums including all $\ell _{p}^{N}$-sums, $1<p\leq \infty $, but not under
$\ell _{1}^{N}$-direct sums (see \cite{La2}). In 1999, B. Sims and M. Smyth \cite{SiSm} proved that
both property (P) and asymptotic (P) are preserved under finite direct sums with monotone norms,
see Section 2 for the relevant definitions. Nowadays, there are many results concerning permanence
properties of normal structure and conditions which imply normal structure (see \cite{Do2, SiSm}),
but only few papers treat a general case of permanence of FPP,
see \cite{DKK, Ku, MaPiXu, Wi3} and
references therein.

In Section 3 we prove two quite general fixed point theorems for direct sums. Theorem \ref{Th1}
states that if $X$ has the weak Banach--Saks property and WFPP, and $Y$ has property asymptotic
(P), then $X\oplus Y$, endowed with a strictly monotone norm, has WFPP. This is a strong extension
of the second named author's results \cite{Wi3}. A combination of the arguments contained in the
proof of Theorem \ref{Th1} with the ideas of P.-K. Lin \cite{Li1} enables us to obtain in Theorem
\ref{Th2} the same conclusion if $X$ has a $1$-unconditional basis, see also a remark at the end of
the paper.

\section{Preliminaries}

Let $X$ be a Banach space and $D\subset X$ be a nonempty set. Given $r>0$, we put
$$ B(D,r)=\{x\in X: \| x-y\| \leq r \text{ for some } y\in D\}. $$
If $D=\{x_0\}$ for some $x_0\in X$, then this is just the closed ball $B(x_0,r)$.

The following construction is crucial for many existence fixed point theorems for nonexpansive
mappings. Assume that there exists a nonexpansive mapping $T:C\rightarrow C$ without a fixed point,
where $C$ is a nonempty weakly compact convex subset of a Banach space $X$. Then, by
the Kuratowski-Zorn lemma,
we obtain a convex and weakly compact set $K\subset C$ which is minimal invariant under $T$ and
which is not a singleton. It follows from the Banach contraction principle that $K$ contains an
approximate fixed point sequence $( x_{n})$ for $T$, i.e.,
\begin{equation*}
\lim_{n\rightarrow \infty }\| Tx_{n}-x_{n}\| =0.
\end{equation*}%

The following lemma was proved independently by K. Goebel \cite{Go} and L.~Karlovitz \cite{Ka}.
\begin{lemma}
\label{GoKa}Let $K$ be a minimal invariant set for a nonexpansive mapping $T$. If $(x_{n})$ is an
approximate fixed point sequence for $T$ in $K$, then%
\begin{equation*}
\lim_{n\rightarrow \infty }\| x_{n}-x\| =\diam K
\end{equation*}%
for every $x\in K$.
\end{lemma}

The above lemma can be reformulated in terms of Banach space ultraproducts as follows (see, e.g.,
\cite{AkKh, Si}). Let $\mathcal{U}$ be a free
ultrafilter on $\mathbb{N}$. The ultrapower $\widetilde{X}:=(X)_{\mathcal{U}%
} $ of a Banach space $X$ is the quotient space of
\begin{equation*}
l_{\infty }(X)=\left\{ (x_{n}):x_{n}\in X\text{ for all }n\in \mathbb{N}\text{ and }%
\| (x_{n})\| =\sup_{n}\| x_{n}\| <\infty \right\}
\end{equation*}%
by
\begin{equation*}
\mathcal{N=}\left\{ (x_{n})\in l_{\infty }(X):\lim_{n\to \mathcal{U}}\| x_{n}\| =0\right\} .
\end{equation*}%
Here $\lim_{n\to \mathcal{U}}$ denotes the ultralimit over $\mathcal{U}$. One can prove that the
quotient norm on $\widetilde{X}$ is given by
\begin{equation*}
\| (x_{n})_{\mathcal{U}}\| =\lim_{n\to \mathcal{U}}\| x_{n}\| ,
\end{equation*}%
where $(x_{n})_{\mathcal{U}}$ is the equivalence class of $(x_{n})$. It is also clear that $X$ is
isometric to a subspace of $\widetilde{X}$ by the mapping $x\mapsto (x,x,\dots)_{\mathcal{U}}$. We
shall not distinguish between $x$ and $(x,x,\dots)_{\mathcal{U}}$. Let
\begin{equation*}
\widetilde{K}=\left\{ (x_{n})_{\mathcal{U}}\in \widetilde{X}:x_{n}\in K%
\text{ for all }n\in \mathbb{N}\right\} .
\end{equation*}%
We extend the mapping $T$ to $\widetilde{K}$ by setting $\widetilde{T}((x_{n})_{\mathcal{U}%
})=(Tx_{n})_{\mathcal{U}}$. It is not difficult to see that $\widetilde{T}:%
\widetilde{K}\rightarrow \widetilde{K}$ is a well-defined nonexpansive mapping. Moreover, the set
$\Fix \widetilde{T}$ of fixed points of $\widetilde{T}$ is nonempty and consists of all those
points in $\widetilde{K}$ which are represented by sequences $(x_{n})$ in $K$ for which $\lim_{n\to
\mathcal{U}}\| Tx_{n}-x_{n}\| =0$. It follows from the Goebel--Karlovitz Lemma \ref{GoKa} ($K$ is
minimal invariant) that
\begin{equation}
\| x-\widetilde{y}\| =\diam K  \label{formula1}
\end{equation}%
for every $x\in K$ and $\widetilde{y}\in \Fix \widetilde{T}$. Even more can be said.
\begin{lemma}[see Lin \cite{Li1}]
\label{Li} Let $K$ be a minimal invariant set for a nonexpansive mapping $T$. If
$(\widetilde{u}_{k})$ is an approximate fixed
point sequence for $\widetilde{T}$ in $\widetilde{K}$, then%
\begin{equation*}
\lim_{k\rightarrow \infty }\| \widetilde{u}_{k}-x\| =\diam K
\end{equation*}%
for every $x\in K$.
\end{lemma}

We conclude this section with recalling several properties of a Banach space $X$ which are
sufficient for weak normal structure. Let
\begin{equation*}
N( X) =\inf \left\{ \frac{\diam A}{r( A) }\right\} ,
\end{equation*}%
where the infimum is taken over all bounded convex sets $A\subset X$ with $\diam A>0$. Assuming
that $X$ does not have the Schur property, we put
\begin{equation*}
\text{\it WCS}\,( X) =\inf \left\{ \frac{\diam_{a}( x_{n}) }{r_{a}( x_{n}) }\right\} ,
\end{equation*}%
where the infimum is taken over all sequences $( x_{n}) $ which converge to $0$ weakly but not in
norm, see \cite{By}. Here
$$\diam_{a}( x_{n}) =\lim_{n\rightarrow \infty }\sup_{k,l\geq n}\|
x_{k}-x_{l}\| $$ denotes the asymptotic diameter of $( x_{n}) $ and
\begin{equation*}
r_{a}( x_{n}) =\inf \left\{ \limsup_{n\rightarrow \infty }\| x_{n}-x\| :x\in \overline{\conv}(
x_{n}) _{n=1}^{\infty }\right\}
\end{equation*}%
denotes the asymptotic radius of $( x_{n}) $.

We say that a Banach space $X$ has uniform normal structure if $N(X) >1$ and weak uniform normal
structure (or satisfies Bynum's condition) if $\text{\it WCS}\,( X) >1$. Moreover, $X$ is said to
have property (P) if
\begin{equation*}
\liminf_{n\rightarrow \infty }\| x_{n}\| <\diam ( x_{n})_{n=1}^{\infty }
\end{equation*}%
whenever $( x_{n}) $ converges weakly to $0$ and $\diam ( x_{n})_{n=1}^{\infty } >0$, see
\cite{TaXu}, and $X$ has property asymptotic
(P) if%
\begin{equation*}
\liminf_{n\rightarrow \infty }\| x_{n}\| <\diam_{a}( x_{n})
\end{equation*}%
whenever $( x_{n}) $ converges weakly to $0$ and $\diam_{a}( x_{n}) >0$, see \cite{SiSm0}. It is
known (see, e.g., \cite{SiSm}) that
\newline
$N( X) >1$ $\Rightarrow$ $\text{\it WCS}\,( X) >1$ $\Rightarrow$ asymptotic (P) $\Rightarrow$ (P)
$\Rightarrow$ weak normal structure.

\section{Results}

In the sequel we shall need the following result (see \cite{Do1, SiSm}).
\begin{lemma}
\label{Do}Every bounded sequence $(x_{n})$ in a Banach space $X$ contains a subsequence $( y_{n}) $
such that the following limit exists
\begin{equation*}
\lim_{n,m\rightarrow \infty ,n\neq m}\| y_{n}-y_{m}\| .
\end{equation*}
\end{lemma}

Let us now recall terminology concerning direct sums. A norm $\left\| \cdot \right\|$ on
$\mathbb{R}^{2}$ is said to be monotone if
\begin{equation*}
\| ( x_{1},y_{1}) \|\leq \| ( x_{2},y_{2}) \|
\end{equation*}%
whenever $0\leq x_{1}\leq x_{2}$, $0\leq y_{1}\leq y_{2}$. A norm $\left\| \cdot \right\|$ is said
to be strictly monotone if
\begin{equation*}
\| ( x_{1},y_{1}) \|<\| ( x_{2},y_{2}) \|
\end{equation*}%
whenever $0\leq x_{1}\leq x_{2}$, $0\leq y_{1}<y_{2}$ or $0\leq
x_{1}<x_{2}$, $0\leq y_{1}\leq y_{2}$. It is easy to see that $\ell _{p}^{2}$%
-norms, $1\leq p<\infty ,$ are strictly monotone.

We shall tacitly assume that
\begin{equation}\label{norm}
\|(1,0)\| =1= \|(0,1)\|.
\end{equation}
This does not result in loss of generality because given a strictly monotone norm
$\left\|\cdot\right\|$ on $\mathbb{R}^2$, we can find another strictly monotone norm
$\left\|\cdot\right\|_1$ such that the spaces $(\mathbb{R}^2,\left\|\cdot\right\|_1)$,
$(\mathbb{R}^2,\left\|\cdot\right\|)$ are isometric and $\left\|\cdot\right\|_1$ satisfies
(\ref{norm}). Moreover, all conditions appearing in our results are isometric invariant.

Let $Z$ be a normed space $(\mathbb{R}^{2},\left\| \cdot \right\|_Z)$. We shall write
$X\oplus_{Z}Y$ for the $Z$-direct sum of Banach spaces $X,Y$ with the norm $\| ( x,y) \| =\| ( \|
x\| ,\| y\|) \|_Z$, where $( x,y) \in X\times Y$.

\begin{lemma}\label{y=0}
Let $X\oplus_{Z}Y$ be a direct sum of Banach spaces $X$, $Y$ with respect to a strictly monotone
norm. Assume that $Y$ has property asymptotic (P), the vectors $v_n=(x_n,y_n)\in X\oplus_{Z}Y$ tend
weakly to 0 and
$$\lim_{n,m\to \infty ,n\neq m}\|v_n-v_m\|=\lim_{n\to \infty}\|v_n\|.$$
Then $\lim_{n\to \infty}\|y_n\|=0$.
\end{lemma}
\begin{proof}
Suppose that the sequence $(y_n)$ does not converge to 0. Then we can assume that the following
limits exist
$$\lim_{n\to \infty }\| x_{n}\|,\quad \lim_{n,m\to \infty ,n\neq m}\| x_{n}-x_{m}\|,\quad
\lim_{n\to \infty }\| y_{n}\|,\quad \lim_{n,m\to \infty ,n\neq m}\| y_{n}-y_{m}\|
$$
(see Lemma \ref{Do}) and $\lim_{n\to \infty }\| y_{n}\|>0$. The sequence $(x_n)$ and $(y_n)$
converges weakly to 0 in $X$ and $Y$, respectively. It follows that
\begin{align*}
\lim_{n\to \infty }\| x_{n}\|&\leq \lim_{n,m\to \infty ,n\neq m}\| x_{n}-x_{m}\|,\\
\lim_{n\to \infty }\| y_{n}\|&\leq \lim_{n,m\to \infty ,n\neq m}\| y_{n}-y_{m}\|.
\end{align*}
Hence
\begin{align*}
\lim_{n\to \infty }\|v_{n} \|&=\|(\lim_{n\to \infty }\| x_{n}\|,\lim_{n\to \infty }\| y_{n}\|) \|_Z
\\
&\leq \|(\lim_{n,m\to \infty ,n\neq m}\| x_{n}-x_{m}\| ,\lim_{n,m\to \infty ,n\neq m}\|
y_{n}-y_{m}\| ) \|_Z
\\ &=\lim_{n,m\to \infty ,n\neq m}\|v_{n} -v_{m}\|=\lim_{n\to \infty }\|v_{n} \|.
\end{align*}
The norm $\left\|\cdot\right\|_Z$ is strictly monotone, so
\begin{equation*}
\lim_{n\to \infty }\| y_{n}\|=\lim_{n,m\to \infty ,n\neq m}\| y_{n}-y_{m}\|= \diam_a ( y_{n})
\end{equation*}
which contradicts our assumption that $Y$ has property asymptotic (P).
\end{proof}

Let $X$ be a Banach space and $x,y\in X$. By the metric segment with the endpoints $x,y$ we mean
the set
$$ S(x,y)=\{z\in X: \|x-z\|+\|z-y\|=\|x-y\|\}. $$
Clearly, $S(x,y)$ contains the algebraic segment $\conv\{x,y\}$.
\begin{lemma}\label{ms}
Let $X\oplus_{Z}Y$ be a direct sum of Banach spaces $X$, $Y$ with respect to a strictly monotone
norm and $\mathcal{U}$ be a free ultrafilter on $\mathbb{N}$. Let $X_0$ denote the set of all
elements of $(X\oplus_{Z}Y)_\mathcal{U}$ of the form $((x_n,0))_\mathcal{U}$ where $(x_n)\in
\ell_\infty(X)$. If $u,v\in X_0$, then $S(u,v)\subset X_0$.
\end{lemma}

\begin{proof}
Let $u=((x_n,0))_\mathcal{U}$, $v=((y_n,0))_\mathcal{U}$ and $z=((a_n,b_n))_\mathcal{U}$ where
$(x_n)$, $(y_n)$, $(a_n)\in \ell_\infty(X)$, $(b_n)\in \ell_\infty(Y)$. Assume that $\lim_{n\to
\mathcal{U}}\|b_n\|>0$ and $z\in S(u,v)$. Since the norm $\left\|\cdot\right\|_Z$ is strictly
monotone,
\begin{align*}
\lim_{n\to \mathcal{U}}\|x_n&-y_n\|=\|u-v\|=\|u-z\|+\|z-v\|\\
&=\|(\lim_{n\to \mathcal{U}}\|x_n-a_n\|, \lim_{n\to \mathcal{U}}\|b_n\|)\|_Z+
\|(\lim_{n\to \mathcal{U}}\|y_n-a_n\|, \lim_{n\to \mathcal{U}}\|b_n\|)\|_Z\\
&>\|(\lim_{n\to \mathcal{U}}\|x_n-a_n\|, 0)\|_Z+ \|(\lim_{n\to \mathcal{U}}\|y_n-a_n\|, 0)\|_Z\\
&=\lim_{n\to \mathcal{U}}\|x_n-a_n\|+\lim_{n\to \mathcal{U}}\|a_n-y_n\|\geq\lim_{n\to
\mathcal{U}}\|x_n-y_n\|.
\end{align*}
This contradiction shows that $\lim_{n\to \mathcal{U}}\|b_n\|=0$ and consequently,
$z=((a_n,0))_\mathcal{U}$.
\end{proof}

Recall that a Banach space $X$ is said to have the weak Banach--Saks property if each weakly null
sequence $( w_{n}) $ in $X$ admits a subsequence $( x_{n}) $ whose arithmetic means converge to $0$
in norm, i.e.,
$$
\lim_{n\rightarrow \infty }\left\| \frac{1}{n}\sum_{k=1}^n x_{k}\right\| =0 .
$$
S. A. Rakov \cite{Ra} proved a result which can be formulated in the following way
(see also \cite{ErMa, FiSu}).
If $(w_n)$ is a
weakly null sequence in a Banach space $X$ with the weak Banach--Saks property, then there is a
subsequence $( x_{n}) $ of $(w_n)$ such that
\begin{equation}\label{ces0}
\lim_{m\to \infty} \sup\left\{\left\|\frac{1}{m}\sum_{i=1}^m x_{p_i}\right\|:
p_1<p_2<\dots<p_m\right\}=0.
\end{equation}

In the proof of the next theorem we shall use the following well-known construction. Let $C$ be a
nonempty convex closed subset of a Banach space $X$ and consider a continuous mapping $T:C\to C$.
Given a separable subset $D$ of $C$, we set $C_1=\conv D$ and $C_{n+1}=\conv(C_n\cup T(C_n))$ for
$n\in \mathbb{N}$. It is easy to see that the set
$$ C(D)=\overline{\bigcup_{n\in \mathbb{N}} C_n} $$
is closed, convex, separable and $T$-invariant. Actually, $C(D)$ is the smallest closed convex
$T$-invariant set containing $D$.

\begin{theorem}
\label{Th1}Let $X$ be a Banach space with the weak Banach--Saks property and WFPP. If $Y$ has
property asymptotic (P), then $X\oplus _{Z}Y$, endowed with a strictly monotone norm, has WFPP.
\end{theorem}

\begin{proof}
Assume that $X\oplus _{Z}Y$ does not have WFPP. Then, there exists a weakly compact convex subset
$C$ of $X\oplus _{Z}Y$ and a nonexpansive mapping $T:C\rightarrow C$ without a fixed point. By the
standard argument described in Section 2, there exists a convex and weakly compact set $K\subset C$
which is minimal invariant under $T$. Let $( w_n ) $ be an approximate fixed point sequence for $T$
in $K$. Without loss of generality we can assume that $\diam K=1$ and $( w_n) $ converges weakly to
$(0,0) \in K$. In view of Lemma \ref{Do} we can assume that the double limit $\lim_{n,m\to \infty
,n\neq m}\|w_n-w_m\|$ exists. From Lemma \ref{GoKa} it follows that
\begin{equation}\label{dseq}
\lim_{n,m\to \infty ,n\neq m}\|w_n-w_m\|=1= \lim_{n\to \infty}\|w_n\|.
\end{equation}
Since $X$ has the weak Banach--Saks property, we can find a subsequence $(x_n)$ of $(w_n)$ for
which condition (\ref{ces0}) holds.

We shall construct by induction a sequence $(n_k^1), (n_k^2),\dots$ of increasing sequences of
natural numbers and an ascending sequence $(D_n)$ of subsets of $\widetilde{K}$ such that for every
$m\in \mathbb{N}$ the following conditions hold
\begin{enumerate}
\item[(i)] the set $\mathbb{N}\setminus \bigcup_{i=1}^{m}A_i$ is infinite and contains $A_{m+1}$
where $A_i=\{n_k^i: k\in \mathbb{N}\}$, \item[(ii)] $D_1=\{v_1\}$ and $D_{m+1} \subset
\overline{\bigcup_{y\in D_m}S(y,v_{m+1})}$, \item[(iii)] $D_m$ is closed, convex, separable,
$\widetilde{T}$-invariant and $\conv\{v_1,\dots,v_{m}\}\subset D_{m}$
\end{enumerate}
where $v_i=(x_{n_k^i})_\mathcal{U}$ for every $i\in \mathbb{N}$.

To this end we put $n_k^1=2k-1$ for every $k\in \mathbb{N}$. Suppose now that we have desirable
sequences $(n_k^1),\dots, (n_k^m)$ and sets $D_1,\dots, D_m$. Then the set $A=\mathbb{N}\setminus
\bigcup_{i=1}^{m}A_i$ is infinite. Let $\{u_n:n\in \mathbb{N}\}$ be a dense subset of $D_m$. We
have $u_k=(y_n^k)_\mathcal{U}$ for some sequence $(y_n^k)$ in $K$. Using Lemma~\ref{GoKa}, for
every $k\in \mathbb{N}$ we find $n_k^{m+1}\in A$ so that $\|y_k^i-x_{n_k^{m+1}}\|\geq
1-\frac{1}{2^k}$ for $i=1,\dots,k$, the sequence $(n_k^{m+1})$ is increasing and the set
$A\setminus A_{m+1}$ is infinite. We have
$$\|u_i-v_{m+1}\|=\lim_{k\to\mathcal{U}}\|y_k^i-x_{n_k^{m+1}}\|=1$$
for every $i\in \mathbb{N}$. It clearly follows that $\|u-v_{m+1}\|=1$ for every $u\in D_m$.

We put $D_{m+1}= C(D_m\cup \{v_{m+1}\})$. To show that (ii) is satisfied observe that the set
$E=\bigcup_{y\in D_m}S(y,v_{m+1})$ is convex and $\widetilde{T}$-invariant. Indeed, if $u_1,u_2\in
E$, then there are $y_1,y_2\in D_m$ such that
$$ \|y_i-u_i\| +\|u_i-v_{m+1}\| =\|y_i-v_{m+1}\|=1$$
for $i=1,2$. Given $t\in [0,1]$, we have $(1-t)y_1+ty_2\in D_m$ and therefore,
\begin{align*}
&\|(1-t)y_1+ty_2-((1-t)u_1+tu_2)\| +\|(1-t)u_1+tu_2-v_{m+1}\|\\&\leq
(1-t)(\|y_1-u_1\|+\|y_1-v_{m+1}\|) +t(\|y_2-u_2\|+\|y_2-v_{m+1}\|)\\
&=(1-t)\|y_1-v_{m+1}\|+t\|y_2-v_{m+1}\|=1\\
&=\|(1-t)y_1+ty_2-v_{m+1}\|.
\end{align*}
This shows that $(1-t)u_1+tu_2\in E$. Moreover, $\widetilde{T}v_{m+1}=v_{m+1}$ and
$\widetilde{T}y_1\in D_m$, so
\begin{align*}
\|\widetilde{T}y_1-\widetilde{T}u_1\| +\|\widetilde{T}u_1-v_{m+1}\| &\leq\|y_1-u_1\|
+\|u_1-v_{m+1}\|=1\\&=\|\widetilde{T}y_1-v_{m+1}\|.
\end{align*}
Therefore, $\widetilde{T}u_1\in E$. Consequently, $E$ is convex, $\widetilde{T}$-invariant and
$D_m\cup \{v_{m+1}\}\subset E$ which easily gives us condition (ii). Condition (iii) is obvious.

We put $D=\overline{\bigcup_{m\in \mathbb{N}}D_m}$. Then $\frac{1}{m}\left(\sum_{i=1}^m v_i\right)
\in D$ for every $m\in \mathbb{N}$ and from (\ref{ces0}) we see that
$$ \lim_{m\to \infty}\left\|\frac{1}{m}\left(\sum_{i=1}^m v_i\right)\right\| =\lim_{m\to \infty}
\lim_{k\to\mathcal{U}} \left\|\frac{1}{m}\left(\sum_{i=1}^m x_{n_k^i}\right)\right\|=0.$$ This
shows that $(0,0)\in D$ and consequently $M=D\cap K\neq \emptyset$. Clearly, $M$ is closed, convex
and $\widetilde{T}$-invariant.

Let $X_0$ denote the set of all elements of $(X\oplus_{Z}Y)_\mathcal{U}$ of the form
$((z_n,0))_\mathcal{U}$ where $(z_n)\in \ell_\infty(X)$. In view of (\ref{dseq}) and Lemma
\ref{y=0}, $v_n\in X_0$ for every $n\in \mathbb{N}$. Using (ii) and Lemma \ref{ms}, one can now
easily show that $D_n\subset X_0$ for every $n\in \mathbb{N}$. Hence $D\subset X_0$ and
consequently $M\subset X_0$. We can therefore identify $M$ with a subset of $X$. Since $X$ has
WFPP, $T$ has a fixed point in $M$ which contradicts our assumption.
\end{proof}

\begin{remark}
The construction of the set $D$ is partly inspired by the arguments in the corrigendum to \cite{Wi3}.
The idea of using metric segments to obtain a $T$-invariant set appeared earlier in \cite{BaNg}.
\end{remark}

It is well known that all superreflexive spaces, $c_{0}$, $\ell_{1}$ as well as $L_{1}[0,1]$ have
the weak Banach--Saks property (see, e.g., \cite{Di}). In metric fixed point theory, the following
coefficient introduced by J. Garc\'{\i}a-Falset \cite{Ga1} plays an important role. Given a Banach
space $X$, we put
$$
R ( X ) =\sup \left\{ \liminf_{n\rightarrow \infty }\| x_{1}+x_{n}\| \right\} ,
$$
where the supremum is taken over all weakly null sequences $( x_{n} )$ in the unit ball $B_{X}$. If
$R ( X ) <2$, then $X$ has the weak Banach--Saks property (see \cite{Ga1}) and WFPP (\cite{Ga2},
see also \cite{Pr2}). For more details about the Banach--Saks property see also \cite{Be, Ca, KPS,
MPS} and references therein.

\begin{corollary}
Let $X$ be a Banach space with $R ( X ) <2$ and $Y$ have weak uniform normal structure. Then
$X\oplus _{Z}Y$, endowed with a strictly monotone norm, has WFPP.
\end{corollary}

It has recently been proved in \cite{GLM2} (see also \cite{MN}) that all uniformly nonsquare Banach
spaces have FPP. Also, uniformly noncreasy spaces introduced in \cite{Pr} are superreflexive and
have FPP. Other examples of superreflexive Banach spaces without normal structure but with FPP are
given by the results in \cite{Wi2, Wi3}.

\begin{corollary}
Let $X$ be a uniformly nonsquare or uniformly noncreasy Banach space and let $Y$ have weak uniform
normal structure. Then $X\oplus _{Z}Y$, endowed with a strictly monotone norm,
has WFPP.
\end{corollary}

In our next result we deal with Banach spaces admitting 1-unconditional bases. Recall that a
Schauder basis $( e_{n})$ of a Banach space $X$ is said to be an unconditional basis provided that
for every choice of signs $( \epsilon _{n}) $, $\epsilon _{n}=\pm 1,$ the series $\sum \epsilon
_{n}\alpha _{n}e_{n}$ converges whenever $\sum \alpha _{n}e_{n}$ converges. Then the supremum
\begin{equation*}
\lambda =\sup \left\{ \left\| \sum_{n=1}^{\infty }\epsilon _{n}\alpha _{n}e_{n}\right\| :\left\|
\sum_{n=1}^{\infty }\alpha _{n}e_{n}\right\| =1,\ \epsilon _{n}=\pm 1\right\}
\end{equation*}%
is finite and it is called the unconditional constant of $( e_{n}) $ (see \cite[p.~18]{LT}). In
this case we say that the basis $( e_{n}) $ is $\lambda$-unconditional. Given a nonempty set
$F\subset \mathbb{N}$, we put
$$ P_Fx= \sum_{n\in F}\alpha _{n}e_{n}$$
where $x= \sum_{n=1}^{\infty }\alpha _{n}e_{n}$. Clearly, $P_F$ is a linear projection and $\| P_F
x\|\leq \lambda \| x\|$ for every $x\in X$.

It is well known that $c_{0}$, $\ell _{p},1\leq p<\infty $ have 1-unconditional bases and the same
is true for the space $X_{\beta }$, which is $\ell _{2}$ endowed with the norm
\begin{equation*}
\| x\| _{\beta }=\max \left\{ \| x\| _{2},\beta \| x\| _{\infty }\right\}.
\end{equation*}%

If we combine the arguments from the first part of the proof of Theorem~\ref{Th1} with the ideas of
Lin \cite{Li1} we obtain the following result.

\begin{theorem}
\label{Th2}Let $X$ be a Banach space with a 1-unconditional basis and let $ Y $ have property
asymptotic (P). Then $X\oplus _{Z}Y$, endowed with a strictly monotone norm, has WFPP.
\end{theorem}

\begin{proof}
Assume that $X\oplus _{Z}Y$ does not have WFPP. Then, there exists a weakly compact convex subset
$K$ of $X\oplus _{Z}Y$ which is minimal invariant under a nonexpansive mapping $T$. Arguing as in
the proof of Theorem~\ref{Th1}, we can assume that an approximate fixed point sequence
$( ( x_{n},y_{n}) ) $ for $T$ converges weakly to $%
( 0,0) \in K$ and that $( y_{n}) $ converges strongly to $0$. Passing to a subsequence, we can
therefore assume that
\begin{equation*}
\lim_{n\rightarrow \infty }\| x_{n}-x_{n+1}\| =\diam K=1.
\end{equation*}

We can now follow the argument from \cite[Theorem 1]{Li1}. Let $(e_n)$ be the 1-unconditional basis
of $X$. By passing to subsequences, we can assume that there exists a sequence $( F_{n}) $ of
intervals of $\mathbb{N}$ such that $\max F_n < \min F_{n+1}$ for every $n\in \mathbb{N}$,
\begin{equation*}
\lim_{n\rightarrow \infty }\left\| P_{F_{n}}x_{n}-x_{n}\right\| =0 \text{\ \ and
}\lim_{n\rightarrow \infty }\left\| P_{F_{n}}x_{n+1}\right\| =\lim_{n\rightarrow \infty }\left\|
P_{F_{n+1}}x_{n}\right\| =0
\end{equation*}%
where $P_{F_k}$ are the projections associated to the basis $(e_n)$. Clearly, $P_{F_{n}}\circ
P_{F_{m}}=0$ if $n\neq m$,
$$ \lim_{n\rightarrow \infty }\left\| P_{F_{n}}x_{n}\right\| =\lim_{n\rightarrow \infty }
\| x_{n}\| =1 \text{ and } \lim_{n\rightarrow \infty }\left\| P_{F_{n}}x\right\| = 0$$ for every
$x\in X$.

Let $\mathcal{U}$ be a free ultrafilter on $\mathbb{N}$ and define projections $\widetilde{P}$,
$\widetilde{Q}$ on $\widetilde{X}$ by
\begin{equation*}
\widetilde{P}( u_{n}) _{\mathcal{U}}=\left( P_{F_{n}}u_{n}\right) _{%
\mathcal{U}} ,\ \widetilde{Q}( u_{n}) _{\mathcal{U}}=\left( P_{F_{n+1}}u_{n}\right) _{\mathcal{U}}.
\end{equation*}%
Put $\widetilde{y}=( x_{n}) _{\mathcal{U}}$, $\widetilde{z}=( x_{n+1}) _{\mathcal{U}}$. Then
\begin{equation}
\widetilde{P}\widetilde{y}=\widetilde{y},\ \widetilde{Q}\widetilde{z}=\widetilde{z}\  \text{ and }\
\widetilde{P}\widetilde{z}=\widetilde{Q}\widetilde{y}=\widetilde{P}x=\widetilde{Q}x=0
\label{formula2}
\end{equation}%
for every $x\in K$. Since $( e_{n}) $ is 1-unconditional,%
\begin{equation*}
\left\| \widetilde{y}+\widetilde{z}\right\| =\left\| \widetilde{y}-\widetilde{z}%
\right\| =1.
\end{equation*}%
Let $\widetilde{v}_1=(( x_{n},0) )_{\mathcal{U}}$, $\widetilde{v}_2= (( x_{n+1},0) )_{\mathcal{U}}$
and
\begin{equation*}
D=B\bigl( \widetilde{v}_1,\tfrac{1}{2}\bigr) \cap B\bigl( \widetilde{v}_2,\tfrac{1}{2}\bigr) \cap
B\bigl( K,\tfrac{1}{2}\bigr) \cap \widetilde{K}.
\end{equation*}

Note that $D\neq \emptyset $ because $\frac{1}{2}\left( \widetilde{v}_1+\widetilde{v}_2\right) \in
D$. Moreover $D$ is closed, convex and $\widetilde{T}( D) \subset D$. Hence, there exists an
approximate fixed point sequence for $\widetilde{T}$ in $D$. Fix an element
$((u_{n},w_{n}))_{\mathcal{U}}\in D$. The set $D$ is contained in the metric segment
$S\left(\widetilde{v}_1, \widetilde{v}_2\right)$. Lemma \ref{ms} shows therefore that
$((u_{n},w_{n})) _{\mathcal{U}}= ( (u_{n},0))_{\mathcal{U}}$. Moreover, from the definition of $D$,
there exists $( x,y) \in K$ such that
\begin{equation*}
\left\| \widetilde{u}-x\right\| \leq \left\| ( (u_{n},0))_{\mathcal{U}}-( x,y) \right\| \leq
\frac{1}{2}
\end{equation*}%
where $\widetilde{u}=( u_{n}) _{\mathcal{U}}$. Hence, with use of (\ref{formula2}), we obtain
\begin{align*}
&\left\| ( (u_{n},w_{n}))_{\mathcal{U}}\right\| =\left\| ( (u_{n},0))_{\mathcal{U}}\right\|
=\left\| \widetilde{u}\right\| \\ &=\frac{1}{2}\left\| \left( \widetilde{u}-\widetilde{P}\widetilde{u}%
\right) +\left( \widetilde{u}-\widetilde{Q}\widetilde{u}\right) +\left( \widetilde{P}\widetilde{u%
}+\widetilde{Q}\widetilde{u}\right) \right\|
\\
&\leq \frac{1}{2}\left(\left\| \Bigl( I-\widetilde{P}\Bigr) ( \widetilde{u}-\widetilde{y}) \right\|
+\left\| \Bigl( I-\widetilde{Q}\Bigr) ( \widetilde{u}- \widetilde{z}) \right\| +\left\| \Bigl(
\widetilde{P}+\widetilde{Q}\Bigr) ( \widetilde{u}-x) \right\|\right)
\\
&\leq \frac{1}{2}\left( \left\|  \widetilde{u}-\widetilde{ y} \right\| +\left\| \widetilde{u}-
\widetilde{z} \right\| +\left\| \widetilde{u}-x \right\|\right) \leq\frac{3}{4},
\end{align*}%
which contradicts Lemma \ref{Li}.
\end{proof}

\begin{remark}
\label{Rem1}In fact, just as in the proof of \cite[Theorem 2]{Li1}, the argument works if $X$ has
an unconditional basis with the unconditional constant $\lambda <\bigl(\sqrt{33}-3\bigr) /2$. Also,
we can adopt the reasoning of M. A. Khamsi \cite{Kh1} (see also \cite[Theorem 4.1]{AkKh}) to obtain
the same conclusion if $X$ is the James quasi-reflexive space.
\end{remark}

Let us recall that there is a separable uniformly convex space which does not embed into a space with an
unconditional basis (see \cite{MR}) and there is a Banach space with a 1-unconditional basis which
does not have the weak Banach--Saks property (see \cite{Ber, CaGo}). This shows that Theorems \ref{Th1}
and \ref{Th2} are entirely independent of each other.

\end{document}